\documentclass[12pt]{article}

\usepackage{amsmath,amssymb,amsthm,amscd,a4wide}
\usepackage{graphicx}

\begin{document}

\newcommand{\End}{{\rm{End}\ts}}
\newcommand{\Hom}{{\rm{Hom}}}
\newcommand{\Mat}{{\rm{Mat}}}
\newcommand{\ad}{{\rm{ad}\ts}}
\newcommand{\ch}{{\rm{ch}\ts}}
\newcommand{\chara}{{\rm{char}\ts}}
\newcommand{\diag}{ {\rm diag}}
\newcommand{\pr}{^{\tss\prime}}
\newcommand{\non}{\nonumber}
\newcommand{\wt}{\widetilde}
\newcommand{\wh}{\widehat}
\newcommand{\ot}{\otimes}
\newcommand{\la}{\lambda}
\newcommand{\ls}{\ts\lambda\ts}
\newcommand{\La}{\Lambda}
\newcommand{\De}{\Delta}
\newcommand{\al}{\alpha}
\newcommand{\be}{\beta}
\newcommand{\ga}{\gamma}
\newcommand{\Ga}{\Gamma}
\newcommand{\ep}{\epsilon}
\newcommand{\ka}{\kappa}
\newcommand{\vk}{\varkappa}
\newcommand{\vt}{\vartheta}
\newcommand{\si}{\sigma}
\newcommand{\vp}{\varphi}
\newcommand{\de}{\delta}
\newcommand{\ze}{\zeta}
\newcommand{\om}{\omega}
\newcommand{\ee}{\epsilon^{}}
\newcommand{\su}{s^{}}
\newcommand{\hra}{\hookrightarrow}
\newcommand{\ve}{\varepsilon}
\newcommand{\ts}{\,}
\newcommand{\vac}{\mathbf{1}}
\newcommand{\di}{\partial}
\newcommand{\qin}{q^{-1}}
\newcommand{\tss}{\hspace{1pt}}
\newcommand{\Sr}{ {\rm S}}
\newcommand{\U}{ {\rm U}}
\newcommand{\BL}{ {\overline L}}
\newcommand{\BE}{ {\overline E}}
\newcommand{\BP}{ {\overline P}}
\newcommand{\AAb}{\mathbb{A}\tss}
\newcommand{\CC}{\mathbb{C}\tss}
\newcommand{\KK}{\mathbb{K}\tss}
\newcommand{\QQ}{\mathbb{Q}\tss}
\newcommand{\SSb}{\mathbb{S}\tss}
\newcommand{\ZZ}{\mathbb{Z}\tss}
\newcommand{\X}{ {\rm X}}
\newcommand{\Y}{ {\rm Y}}
\newcommand{\Z}{{\rm Z}}
\newcommand{\Ac}{\mathcal{A}}
\newcommand{\Lc}{\mathcal{L}}
\newcommand{\Mc}{\mathcal{M}}
\newcommand{\Pc}{\mathcal{P}}
\newcommand{\Qc}{\mathcal{Q}}
\newcommand{\Tc}{\mathcal{T}}
\newcommand{\Sc}{\mathcal{S}}
\newcommand{\Bc}{\mathcal{B}}
\newcommand{\Ec}{\mathcal{E}}
\newcommand{\Fc}{\mathcal{F}}
\newcommand{\Hc}{\mathcal{H}}
\newcommand{\Uc}{\mathcal{U}}
\newcommand{\Vc}{\mathcal{V}}
\newcommand{\Wc}{\mathcal{W}}
\newcommand{\Yc}{\mathcal{Y}}
\newcommand{\Ar}{{\rm A}}
\newcommand{\Br}{{\rm B}}
\newcommand{\Ir}{{\rm I}}
\newcommand{\Fr}{{\rm F}}
\newcommand{\Jr}{{\rm J}}
\newcommand{\Mr}{{\rm M}}
\newcommand{\Or}{{\rm O}}
\newcommand{\GL}{{\rm GL}}
\newcommand{\SL}{{\rm SL}}
\newcommand{\Spr}{{\rm Sp}}
\newcommand{\Rr}{{\rm R}}
\newcommand{\Zr}{{\rm Z}}
\newcommand{\gl}{\mathfrak{gl}}
\newcommand{\middd}{{\rm mid}}
\newcommand{\ev}{{\rm ev}}
\newcommand{\Pf}{{\rm Pf}}
\newcommand{\Norm}{{\rm Norm\tss}}
\newcommand{\oa}{\mathfrak{o}}
\newcommand{\spa}{\mathfrak{sp}}
\newcommand{\osp}{\mathfrak{osp}}
\newcommand{\bgot}{\mathfrak{b}}
\newcommand{\kgot}{\mathfrak{k}}
\newcommand{\g}{\mathfrak{g}}
\newcommand{\h}{\mathfrak h}
\newcommand{\n}{\mathfrak n}
\newcommand{\z}{\mathfrak{z}}
\newcommand{\Zgot}{\mathfrak{Z}}
\newcommand{\p}{\mathfrak{p}}
\newcommand{\sll}{\mathfrak{sl}}
\newcommand{\psl}{\mathfrak{psl}}
\newcommand{\agot}{\mathfrak{a}}
\newcommand{\qdet}{ {\rm qdet}\ts}
\newcommand{\Ber}{ {\rm Ber}\ts}
\newcommand{\HC}{ {\mathcal HC}}
\newcommand{\cdet}{ {\rm cdet}}
\newcommand{\tr}{ {\rm tr}}
\newcommand{\gr}{ {\rm gr}}
\newcommand{\str}{ {\rm str}}
\newcommand{\loc}{{\rm loc}}
\newcommand{\Gr}{{\rm G}}
\newcommand{\sgn}{ {\rm sgn}\ts}
\newcommand{\ba}{\bar{a}}
\newcommand{\bb}{\bar{b}}
\newcommand{\bi}{\bar{\imath}}
\newcommand{\bj}{\bar{\jmath}}
\newcommand{\bk}{\bar{k}}
\newcommand{\bl}{\bar{l}}
\newcommand{\hb}{\mathbf{h}}
\newcommand{\Sym}{\mathfrak S}
\newcommand{\fand}{\quad\text{and}\quad}
\newcommand{\Fand}{\qquad\text{and}\qquad}
\newcommand{\For}{\qquad\text{or}\qquad}
\newcommand{\OR}{\qquad\text{or}\qquad}

\renewcommand{\theequation}{\arabic{section}.\arabic{equation}}

\newtheorem{thm}{Theorem}[section]
\newtheorem{lem}[thm]{Lemma}
\newtheorem{prop}[thm]{Proposition}
\newtheorem{cor}[thm]{Corollary}
\newtheorem{conj}[thm]{Conjecture}
\newtheorem*{mthm}{Main Theorem}
\newtheorem*{mthma}{Theorem A}
\newtheorem*{mthmb}{Theorem B}

\theoremstyle{definition}
\newtheorem{defin}[thm]{Definition}

\theoremstyle{remark}
\newtheorem{remark}[thm]{Remark}
\newtheorem{example}[thm]{Example}

\newcommand{\bth}{\begin{thm}}
\renewcommand{\eth}{\end{thm}}
\newcommand{\bpr}{\begin{prop}}
\newcommand{\epr}{\end{prop}}
\newcommand{\ble}{\begin{lem}}
\newcommand{\ele}{\end{lem}}
\newcommand{\bco}{\begin{cor}}
\newcommand{\eco}{\end{cor}}
\newcommand{\bde}{\begin{defin}}
\newcommand{\ede}{\end{defin}}
\newcommand{\bex}{\begin{example}}
\newcommand{\eex}{\end{example}}
\newcommand{\bre}{\begin{remark}}
\newcommand{\ere}{\end{remark}}
\newcommand{\bcj}{\begin{conj}}
\newcommand{\ecj}{\end{conj}}

\newcommand{\bal}{\begin{aligned}}
\newcommand{\eal}{\end{aligned}}
\newcommand{\beq}{\begin{equation}}
\newcommand{\eeq}{\end{equation}}
\newcommand{\ben}{\begin{equation*}}
\newcommand{\een}{\end{equation*}}

\newcommand{\bpf}{\begin{proof}}
\newcommand{\epf}{\end{proof}}

\def\beql#1{\begin{equation}\label{#1}}

\newcommand{\cO}{{\mathcal O}}
\newcommand{\cQ}{{\mathcal Q}}
\newcommand{\sC}{{C\!\!\!\!/\,}}
\newcommand{\fL}{{\mathfrak L}}
\newcommand{\fR}{{\mathfrak R}}
\newcommand{\fQ}{{\mathfrak Q}}
\newcommand{\fS}{{\mathfrak S}}
\newcommand{\fC}{{\mathfrak C}}
\newcommand{\fP}{{\mathfrak P}}
\newcommand{\fK}{{\mathfrak K}}
\newcommand{\fJ}{{\mathfrak J}}
\newcommand{\fB}{{\mathfrak B}}
\newcommand{\fF}{{\mathfrak F}}
\newcommand{\fe}{{\mathfrak e}}
\newcommand{\fh}{{\mathfrak h}}

\newcommand{\sfR}{{\mathfrak R}\!\!\!\!/\,}
\newcommand{\sfL}{{\mathfrak L}\!\!\!/}
\newcommand{\sfC}{{\mathfrak C}\!\!\!/}
\newcommand{\cT}{\mathcal T}
\newcommand{\Tr}{{\rm Tr}}
\newcommand{\bs}[1]{\boldsymbol{#1}}
\newcommand{\alg}[1]{\mathfrak{#1}}
\newcommand{\el}{\nonumber}
\newcommand{\nln}{\nonumber\\}
\newcommand{\eg}{\widetilde{{\mathfrak g}}}

\newtheorem{dfn}{Definition}[section]
\newtheorem{rem}{Remark}[section]

\title{\Large\bf Representations of centrally
extended\\ Lie superalgebra $\psl(2|2)$}

\author{{Takuya Matsumoto\quad and\quad Alexander Molev}}

\date{} 
\maketitle

\vspace{25 mm}

\begin{abstract}
The symmetries provided by
representations of the centrally extended Lie superalgebra $\psl(2|2)$
are known to play an important role in the spin chain models
originated in the planar anti-de Sitter/conformal field theory correspondence
and one-dimensional Hubbard model.
We give a complete description of
finite-dimensional irreducible
representations of this superalgebra thus extending the work of Beisert
which deals with a generic family of representations.
Our description includes a new class of modules with degenerate
eigenvalues of the central elements. Moreover, we construct explicit bases
in all irreducible representations by applying the techniques of
Mickelsson--Zhelobenko algebras.

%

\end{abstract}


\vspace{35 mm}

\noindent
Institute for Theoretical Physics and Spinoza Institute \\
Utrecht University, Leuvenlaan 4, 3854 CE Utrecht, The Netherlands  \\
t.matsumoto@uu.nl

\vspace{7 mm}

\noindent
School of Mathematics and Statistics\newline
University of Sydney,
NSW 2006, Australia\newline
alexander.molev@sydney.edu.au

%

\newpage

\section{Introduction}
\label{sec:int}
\setcounter{equation}{0}

As discovered by Beisert \cite{b:su, B-yang, b:ab}, certain spin chain models originated in
the planar anti-de Sitter/conformal field theory (AdS/CFT)  correspondence admit
{\em hidden symmetries} provided by the action of the Yangian $\Y(\g)$ associated with the
centrally extended Lie superalgebra
\ben
\g=\psl(2|2)\ltimes \CC^3.
\een
This is a semi-direct product of the simple Lie superalgebra $\psl(2|2)$
of type $A(1,1)$
and the abelian Lie algebra $\CC^3$ spanned by elements $C$, $K$ and $P$ which
are central in $\g$.
Due to the results of \cite{ik:ce},
$\psl(2|2)$ is distinguished among the basic classical Lie superalgebras
by the existence of a three-dimensional central extension. It was pointed out in
\cite{d:nc} that this phenomenon originates in some special Serre relations.
A new $R$-matrix associated with the extended Lie superalgebra
$\g$ is found by Yamane~\cite{y:ce}.
Furthermore, $\g$ can be obtained from the Lie superalgebras
of type $D(2,1;\al)$ by a particular limit with respect to the parameter~$\al$.

The Yangian symmetries of the one-dimensional Hubbard model
associated with $\Y(\g)$ were considered in \cite{B-yang};
they extend those provided by the direct sum of two copies
of the Yangian for
$\sll(2)$ previously found in \cite{uk:ys}.
An extensive review of the Yangian
symmetries in the spin chain models can be found in \cite{t:ys}.

These applications motivate
the study of representations of both
the Lie superalgebra $\g$ and its Yangian. In this paper we aim to prove
a classification theorem for finite-dimensional irreducible representation
of $\g$. Generic representations of $\g$ were already described by
Beisert~\cite{b:ab}. As we demonstrate below, beside these generic modules,
the complete classification includes some degenerate representations
which were not considered in \cite{b:ab}.
In more detail, if $L$ is a finite-dimensional irreducible representation of the Lie
superalgebra $\g$, then each of the central elements $C,K$ and $P$ acts in $L$
as multiplication by a scalar. We will let the lower case letters denote
the corresponding scalars,
\ben
C\mapsto c,\qquad K\mapsto k,\qquad
P\mapsto p.
\een
The Lie superalgebra $\psl(2|2)$ is known to admit a family
of automorphisms parameterized by elements of the group $\SL(2)$,
as described in \cite{gqs:tp}.
As pointed out in \cite{b:ab}, by twisting the action of $\g$ in $L$ by such an
automorphism, we obtain
another irreducible representation of $\g$, where the
values $c,k,p$ are transformed by
\begin{align}
\begin{pmatrix}c & -k \\ p & -c \end{pmatrix}\mapsto
\begin{pmatrix}u & v \\ w & z \end{pmatrix}
\begin{pmatrix}c & -k \\ p & -c \end{pmatrix}
\begin{pmatrix}u & v \\ w & z \end{pmatrix}^{-1},
\non
\end{align}
and complex numbers $u,v,w,z$ satisfy $u\tss z-v\tss w=1$.
An appropriate
transformation of this form brings the $2\times 2$ matrix formed by
$c,k,p$ to the Jordan canonical form. In the case, where the canonical form
is a diagonal matrix,
\ben
\begin{pmatrix}c & -k \\ p & -c \end{pmatrix}\mapsto
\begin{pmatrix}d & \phantom{-}0 \\ 0 & -d \end{pmatrix},
\een
the values of $k$ and $p$ under the twisted action
of $\g$ are zero, and so the twisted module becomes an irreducible representation
of the Lie superalgebra $\sll(2|2)$. Such representations are well-studied;
see e.g. \cite{gqs:tp}, \cite{kkp:fd}, \cite{ps:fd} and \cite{zg:uc} for
an explicit construction of basis vectors and formulas
for the action of the generators in this basis. It is essentially this case
which was considered in \cite{b:ab} in relation with the symmetries of the $S$-matrix
for the AdS/CFT correspondence. The only remaining possibility is the case where
the canonical form is the $2\times 2$ Jordan block,
\ben
\begin{pmatrix}c & -k \\ p & -c \end{pmatrix}\mapsto
\begin{pmatrix}0 & 0 \\ 1 & 0 \end{pmatrix}.
\een
This can only happen when $c^2-p\tss k=0$ so that both eigenvalues
of the matrix are zero. Our goal in this paper is to study the
structure of these representations of $\g$. In what follows we consider
the class of finite-dimensional irreducible representations of $\g$
where both central elements $C$ and $K$ act as the zero operators, while
$P$ acts as the identity operator.
Our main result
is a classification theorem
for such representations of $\g$.

\begin{mthm}
A complete list of pairwise non-isomorphic finite-dimensional
irreducible representations of $\g$ where the central elements
act by $C\mapsto 0$, $K\mapsto 0$, $P\mapsto 1$,
consists of
\begin{enumerate}
\item
the Kac modules $K(m,n)$ with $m,n\in\ZZ_+$ and $m\ne n$,\newline
$\dim K(m,n)=16\tss(m+1)(n+1)$,
\item
the modules $S_n$ with $n\in\ZZ_+$,\quad $\dim S_n=8\tss(n+1)(n+2)$.
\end{enumerate}
\end{mthm}

Here the {\em Kac modules} $K(m,n)$ over $\g$ are defined as the induced modules
from finite-dimensional irreducible representations of the Lie algebra
$\sll(2)\oplus \sll(2)$ in the same way as for simple Lie superalgebras, and $S_n$ are
certain submodules of $K(n,n)$; see Sec.~\ref{subsec:kac} for precise
definitions. Comparing this description with the classification theorem
for representations of the Lie superalgebra $\psl(2|2)$ \cite{k:ct}, note that
nontrivial irreducible quotients of the Kac modules over $\psl(2|2)$
have dimension $4\tss n(n+2)+2$ (they are also known as {\em short multiplets}).
Hence the family of $\g$-modules $S_n$
does not have their counterparts within the class of $\psl(2|2)$-modules.

To give a physical interpretation of the conditions on $c,p$ and $k$, note that
in the original spin chain models \cite{b:su, b:ab}, the scalar $c$
corresponds to the {\em energy} of a particle moving on the spin chain, whereas
$p$ and $k$ correspond to its {\em momenta}.
Thus,  the relation $c^2-p\tss k=\text{const}$ is the dispersion relation of the particle,
and the automorphisms of $\psl(2|2)$ provided by elements of $\SL(2)$
are interpreted as the {\em Lorentz symmetry}
which preserves the dispersion relation.
Therefore, the relation $c^2-p\tss k=0$ describes a massless particle on the {\em light-cone}.
Due to the Main Theorem, the multiplets of the particle on the light-cone are {\em shorter} than
long multiplets \cite{b:ab} and {\em longer} than short multiplets.
Thus, the particles on the light-cone are described by {\em middle} multiplets.

Our arguments are based on the theory of Mickelsson--Zhelobenko algebras
\cite{z:it}. We also apply it to construct
bases of all finite-dimensional irreducible representations of $\g$.
Formulas for the action of the generators of $\g$ in such a basis
can also be found in an explicit form.
Furthermore, this description of representations extends to
the case, where the central elements $P$ and $K$ of $\g$ act as the zero operators,
allowing us to essentially reproduce the results of \cite{kkp:fd},
\cite{ps:fd} and \cite{zg:uc} concerning
representations of $\gl(2|2)$ and $\sll(2|2)$.

This paper is organized as follows.
In Sec.~\ref{sec:ce} we review the centrally extended Lie superalgebra $\g$.
In Sec.~\ref{sec:fd} we describe
finite-dimensional irreducible representation of $\g$.
After introducing the Mickelsson--Zhelobenko algebra in Sec.~\ref{subsec:mz}
we construct a basis of the Kac module by the Mickelsson--Zhelobenko generators
and establish its irreducibility properties.
In Sec.~\ref{subsec:ct} the classification theorem is proved.
Explicit action of the generators on the Kac modules is described
in Sec.~\ref{subsec:ec}.
Appendix~\ref{sec:amz} is devoted to relations in
the Mickelsson--Zhelobenko algebra.
In Appendix~\ref{sec:aar} the action of raising
operators is produced; it is used to prove
irreducibility of the Kac modules.
In Appendix~\ref{sec:abk} the $\sll(2)\oplus\sll(2)$-highest weight vectors
of the Kac modules are listed.

\medskip

We acknowledge the financial support of the Australian Research Council.
T.M.\ would like to thank the hospitality of the School of
Mathematics and Statistics at the University of Sydney.
The most of this work has been done during his stay there.
T.M. also thanks Sanefumi Moriyama, Hiroyuki Yamane
and Kentaroh Yoshida for variable discussions.
T.M.\ is supported by the Netherlands Organization for Scientific
Research (NWO) under the VICI grant 680-47-602.
T.M.'s work is also part of the ERC Advanced grant research programme
No.~246974, ``Supersymmetry: a window to non-perturbative physics"
and of the D-ITP consortium, a program of the NWO that is funded by the
Dutch Ministry of Education, Culture and Science (OCW).

\section{Central extension of Lie superalgebra $\alg{psl}(2|2)$}
\label{sec:ce}
\setcounter{equation}{0}

The general linear Lie superalgebra $\gl(2|2)$ over $\CC$
has the standard basis
$E_{ij}$, $1\leqslant i,j\leqslant 4$. The $\ZZ_2$-grading on $\gl(2|2)$
is defined by setting $\deg E_{ij}=\bi+\bj$, where
we use the notation $\bi=0$ for $1\leqslant i\leqslant 2$
and $\bi=1$ for $3\leqslant i\leqslant 4$.
The commutation relations
have the form
\ben
\big[E_{ij},E_{kl}\tss\big]
=\de_{kj}\ts E_{i\tss l}-\de_{i\tss l}\ts
E_{kj}(-1)^{(\bi+\bj)(\bk+\bl)},
\een
where the square brackets denote the super-commutator. Then
$\sll(2|2)$ is the subalgebra of $\gl(2|2)$ spanned by
the elements
\ben
h_1=E_{11}-E_{22},\qquad h_2=E_{22}+E_{33}, \qquad h_3=E_{33}-E_{44}
\een
and by all elements $E_{ij}$ with $i\ne j$.
We have the direct sum
decomposition
\ben
\gl(2|2)=\sll(2|2)\oplus \CC(E_{11}+E_{22}-E_{33}-E_{44}).
\een
Furthermore,
the element
\ben
C=\frac12\ts h_1+h_2-\frac12\ts h_3=\frac12\ts (E_{11}+E_{22}+E_{33}+E_{44})
\een
is central in $\sll(2|2)$, and the simple Lie superalgebra of type $A(1,1)$
is defined as the quotient of $\sll(2|2)$ by the ideal generated by $C$.
This quotient is denoted by $\psl(2|2)$. As in \cite{b:ab}, we will consider
the Lie superalgebra
\ben
\g=\sll(2|2)\ltimes\CC^2 =\psl(2|2)\ltimes\CC^3,
\een
where $\CC^2$ is the abelian Lie algebra with the basis elements
$K$ and $P$, while abelian Lie algebra
$\CC^3$ is spanned by the elements $C,K$ and $P$.
These elements are central in $\g$ and the only nontrivial additional
relations take the form
\begin{align}
&[E_{13}, E_{24}]=-[E_{23}, E_{14}]=K, \label{extk}\\
&[E_{31}, E_{42}]=-[E_{32}, E_{41}]=P. \label{extp}
\end{align}
More precisely, the commutations relations in $\g$ are determined by those
for the basis elements
\beql{gl22}
[E_{ij}, E_{kl}]=\delta_{kj} E_{il}-\delta_{il} E_{kj}(-1)^{(\bi+\bj)(\bk+\bl)}
+\bar \epsilon_{ik}\tss\epsilon_{jl} P+\epsilon_{ik}\tss\bar \epsilon_{jl} K,
\eeq
where the constants $\ep_{ij}$ and $\bar \ep_{ij}$ are zero except for
the values
\ben
\epsilon_{12}=-\epsilon_{21}=1\Fand\bar\epsilon_{34}=-\bar\epsilon_{43}=1.
\een
The Lie subalgebra $\g_0$ of even elements in $\g$ is the direct sum
\ben
\g_0=\kgot\oplus\CC^3,\qquad \kgot=\sll(2)\oplus\sll(2)\subset\psl(2|2),
\een
where the two copies of $\sll(2)$ are spanned by the
elements $E_{12}, E_{21}, h_1$ and $E_{34}, E_{43}, h_3$, respectively.

Given complex numbers $u,v,w,z$ such that $u\tss z-v\tss w=1$,
the corresponding automorphism $\phi:\g\to\g$ mentioned in the Introduction
is determined by the mapping
\beql{autophi}
E_{13}\mapsto u\tss E_{13}+v\tss E_{42},\qquad
E_{42}\mapsto z\tss E_{42}+w\tss E_{13},
\eeq
and the condition that
each element of the subalgebra $\sll(2)\oplus\sll(2)$ is stable under $\phi$.
By \cite{b:ab}, the images of the central elements $C,K,P$
are then found from the matrix relation
\begin{align}
\begin{pmatrix}C & -K \\ P & -C \end{pmatrix}\mapsto
\begin{pmatrix}u & v \\ w & z \end{pmatrix}
\begin{pmatrix}C & -K \\ P & -C \end{pmatrix}
\begin{pmatrix}u & v \\ w & z \end{pmatrix}^{-1}.
\non
\end{align}

\section{Finite-dimensional irreducible representations}
\label{sec:fd}
\setcounter{equation}{0}

As we pointed out in the Introduction, our main focus will be on
finite-dimensional irreducible representations of $\g$,
where the eigenvalues of the central elements are given by
\ben
C\mapsto 0,\qquad K\mapsto 0,\qquad
P\mapsto 1.
\een
This means that we will essentially deal with
the extended Lie superalgebra $\psl(2|2)\oplus\CC P$,
where the only nontrivial additional relations are \eqref{extp}.

From the viewpoint of the spin chain model \cite{b:ab}, these representations
should describe the particle states on the {\em light-cone} since the
{\em dispersion relations}
are given by $c^2-p\tss k=0$, where $c$ and $p,k$ correspond to {\em energy } and {\em momenta}
of the particles, respectively.

\subsection{Mickelsson--Zhelobenko algebras}
\label{subsec:mz}

We will use the {\em Mickelsson--Zhelobenko algebra} $\Z(\g,\kgot)$
associated with the pair $\kgot\subset\g$. An extensive theory
of such algebras was developed in \cite{z:it}; see also
\cite[Ch.~9]{m:yc} and \cite{m:cb} where they were employed
for constructions of bases of Gelfand--Tsetlin type in representations
of classical Lie algebras and superalgebras. To recall the definitions,
denote by $\h$ the Cartan subalgebra of $\kgot$ spanned by the basis
elements $h_1$ and $h_3$. We have
the triangular
decomposition
\ben
\kgot=\kgot^-\oplus\h\oplus\kgot^+,
\een
where
\ben
\kgot^-=\text{span of}\ \{E_{21}, E_{43}\}\Fand
\kgot^+=\text{span of}\ \{E_{12}, E_{34}\}.
\een
Let $\Jr=\U(\g)\ts\kgot^+$ be
the left ideal of $\U(\g)$
generated by $\kgot^+$ and consider the quotient
\ben
\Mr(\g,\kgot)= \U(\g)/\Jr.
\een
The {\em Mickelsson algebra} $\Sr(\g,\kgot)$ is defined by
\ben
\Sr(\g,\kgot)=\{v\in \Mr(\g,\kgot)\ |\ \kgot^+ v=0\}.
\een
Given a finite-dimensional $\g$-module $V$, its subspace
\beql{vplus}
V^+=\{v\in V\ |\ \kgot^+ v=0\}
\eeq
is a $\Sr(\g,\kgot)$-module whose structure largely determines the structure
of $V$; see \cite{z:it} for more details.
Denote by $\Rr(\h)$
the field of fractions of the commutative algebra $\U(\h)$.
The {\em Mickelsson--Zhelobenko algebra} $\Z(\g,\kgot)$ can be defined
as the extension
\beql{extueal}
\Z(\g,\kgot)=\Sr(\g,\kgot)\ot_{\U(\h)} \Rr(\h).
\eeq
As was observed by Zhelobenko (see \cite{z:it}),
the algebraic structure of $\Z(\g,\kgot)$ can be described with
the use of the {\em extremal projector} $p=p(\kgot)$
which is a formal series of elements of $\U(\kgot)$ with
coefficients in $\Rr(\h)$ given by
\begin{multline}
p=\Bigg(1+\sum_{k=1}^\infty \frac{(-1)^k}{k!}E_{21}^kE_{12}^k
\frac{1}{(h_1+2)\cdots(h_1+k+1)}\Bigg) \\[0.3em]
{}\times \Bigg(1+\sum_{l=1}^\infty
\frac{(-1)^l}{l!}E_{43}^lE_{34}^l\frac{1}{(h_3+2)\cdots(h_3+l+1)}\Bigg).
\non
\end{multline}
The operator $p$ has the property $p^2=p$ and satisfies the relations
\beql{eppe}
E_{12}\tss p=p\tss E_{21}=0\Fand E_{34}\tss p=p\tss E_{43}=0.
\eeq
The extremal projector naturally acts on the extension
\ben
\Mr'(\g,\kgot)=\Mr(\g,\kgot)\ot_{\U(\h)} \Rr(\h).
\een
It projects $\Mr'(\g,\kgot)$
onto $\Z(\g,\kgot)$ with the kernel $\kgot^-\Mr'(\g,\kgot)$.
In particular, $\Z(\g,\kgot)=p\ts \Mr'(\g,\kgot)$. Moreover,
the algebra $\Z(\g,\kgot)$ is generated by the elements
\ben
z_{ik}=p\tss E_{ik},\qquad z_{ki}=p\tss E_{ki},\qquad i=1,2\fand k=3,4,
\een
together with
$C,K$ and $P$. We will call the elements $z_{ik}$ and $z_{ki}$
{\em raising {\em and} lowering} operators, respectively.
They are given by the following explicit formulas.

\ble\label{lem:raislower}
The raising operators are found by
\begin{align}
z_{14}&=E_{14}, \el \\
z_{13}&=E_{13}+E_{43}E_{14}\tss\frac{1}{h_3+1}, \el \\
z_{24}&=E_{24}-E_{21}E_{14}\tss\frac{1}{h_1+1}, \el \\
z_{23}&=E_{23}-E_{21}E_{13}\tss\frac{1}{h_1+1}+E_{43}E_{24}\tss\frac{1}{h_3+1}
-E_{21}E_{43}E_{14}\tss\frac{1}{(h_1+1)(h_3+1)},
\el
\end{align}
and the lowering operators are
\begin{align}
z_{41}&=E_{41}+E_{21}E_{42}\tss\frac{1}{h_1+1}-E_{43}E_{31}\tss\frac{1}{h_3+1}
-E_{21}E_{43}E_{32}\tss\frac{1}{(h_1+1)(h_3+1)},   \el \\
z_{31}&=E_{31}+E_{21}E_{32}\tss\frac{1}{h_1+1}, \el \\
z_{42}&=E_{42}-E_{43}E_{32}\tss\frac{1}{h_3+1}, \el \\
z_{32}&=E_{32}.
\el
\end{align}
\ele

\bpf
These expressions follow by the application of the explicit formula for the
extremal projector $p$.
\epf

We will need expressions for
the elements $E_{ik}$ and $E_{ki}$ in terms of the raising and lowering
operators provided by the next lemma.

\ble\label{lem:inve}
We have the relations in $\Mr'(\g,\kgot)$:
\begin{align}
E_{14}&=z_{14}, \el \\
E_{13}&=z_{13}-E_{43}z_{14}\tss\frac{1}{h_3+1}, \el \\
E_{24}&=z_{24}+E_{21}z_{14}\tss\frac{1}{h_1+1}, \el \\
E_{23}&=z_{23}+E_{21}z_{13}\tss\frac{1}{h_1+1}-E_{43}z_{24}\tss\frac{1}{h_3+1}
-E_{21}E_{43}z_{14}\tss\frac{1}{(h_1+1)(h_3+1)},
\el
\end{align}
and
\begin{align}
E_{41}&=z_{41}-E_{21}z_{42}\tss\frac{1}{h_1+1}+E_{43}z_{31}\tss\frac{1}{h_3+1}
-E_{21}E_{43}z_{32}\tss\frac{1}{(h_1+1)(h_3+1)},   \el \\
E_{31}&=z_{31}-E_{21}z_{32}\tss\frac{1}{h_1+1}, \el \\
E_{42}&=z_{42}+E_{43}z_{32}\tss\frac{1}{h_3+1}, \el \\
E_{32}&=z_{32}.
\el
\end{align}
\ele

\bpf
The formulas are immediate from Lemma~\ref{lem:raislower}.
\epf

As follows from \cite{z:it}, the generators of the Mickelsson--Zhelobenko algebra
$\Zr(\g,\kgot)$ satisfy quadratic relations which can be
derived from Lemmas~\ref{lem:raislower} and \ref{lem:inve}.
In particular,
for $i=1,2$ and $k=3,4$ we have
\ben
z_{ik}^2=0\Fand z_{ki}^2=0.
\een
Complete sets of relations in $\Zr(\g,\kgot)$ are listed in Appendix~\ref{sec:amz}.

\subsection{Kac modules}
\label{subsec:kac}

For nonnegative integers $m$ and $n$
we will denote by $L^0(m,n)$ the finite-dimensional irreducible
representation of the Lie algebra $\kgot=\sll(2)\oplus\sll(2)$
with the highest weight $[m,n]$. This representation is generated
by the highest vector $w$ such that
\ben
E_{12}\tss w=E_{34}\tss w=0,\qquad h_1\tss w=m\tss w,\qquad h_3\tss w=n\tss w.
\een
The vectors
\beql{basl}
E_{21}^k E_{43}^l\tss w,\qquad k=0,1,\dots,m,\qquad l=0,1,\dots,n,
\eeq
form a basis of $L^0(m,n)$.
We extend $L^0(m,n)$ to a representation of the subalgebra $\bgot\subset\g$,
spanned by $\g_0=\kgot\oplus\CC^3$ and the elements $E_{ik}$ with $i=1,2$ and $k=3,4$.
These additional elements act as the zero operators, while
$C\mapsto 0$, $K\mapsto 0$ and $P\mapsto 1$.
The corresponding {\em Kac module} $K(m,n)$ is defined
as the induced representation
\beql{kacmodule}
K(m,n)=\U(\g)\ot_{\U(\bgot)} L^0(m,n).
\eeq
Its basis is formed by the vectors
\ben
E_{41}^{\theta_1}\tss E_{31}^{\theta_2}\tss E_{42}^{\theta_3}\tss  E_{32}^{\theta_4}
\tss E_{21}^k E_{43}^l\tss w
\een
where each
$\theta_i$ takes values in $\{0,1\}$ and $k,l$ are as in \eqref{basl}.
In particular,
\ben
\dim K(m,n)=16\tss(m+1)(n+1).
\een

As in \eqref{vplus}
we will write $K^+(m,n)$
for the subspace of $\kgot^+$-invariants in $K(m,n)$.
Note that the action
of the elements $z_{ik}$ and $z_{ki}$ of the Mickelsson--Zhelobenko algebra
$\Z(\g,\kgot)$ in $K^+(m,n)$ is well-defined. The denominators
of these rational functions do not vanish when $h_1$ and $h_3$ are replaced
by the corresponding eigenvalues of weight vectors.
In accordance with \cite{z:it},
a basis of $K^+(m,n)$
can be constructed with the use of ordered products of the lowering operators
$z_{ki}$ with $i=1,2$ and $k=3,4$. Below we use this approach
to write down explicit basis elements.

Suppose first that $m,n\geqslant 2$.
Consider the elements
\beql{baskac}
z_{41}^{\theta_1}\ts
z_{31}^{\theta_2}\ts z_{42}^{\theta_3}\ts z_{32}^{\theta_4}\ts w\in K^+(m,n)
\eeq
with each $\theta_i$ taking values in $\{0,1\}$.
Each element \eqref{baskac} can be interpreted as a path in the
following labelled oriented graph

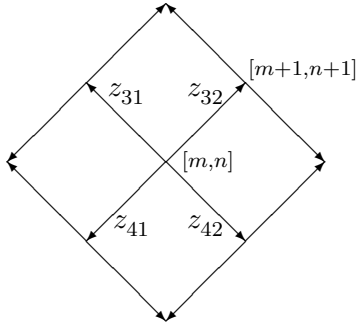
\begin{center}
\begin{picture}(200,165)
\thinlines

\put(100,100){\vector(1,1){30}}
\put(100,100){\vector(1,-1){30}}
\put(100,100){\vector(-1,1){30}}
\put(100,100){\vector(-1,-1){30}}
\put(130,70){\vector(1,1){30}}
\put(130,70){\vector(-1,-1){30}}
\put(130,130){\vector(-1,1){30}}
\put(130,130){\vector(1,-1){30}}
\put(70,130){\vector(1,1){30}}
\put(70,130){\vector(-1,-1){30}}
\put(70,70){\vector(-1,1){30}}
\put(70,70){\vector(1,-1){30}}

\put(108,125){\small $z^{}_{32}$}
\put(78,125){\small $z^{}_{31}$}
\put(108,75){\small $z^{}_{42}$}
\put(80,75){\small $z^{}_{41}$}

\put(106,98){$\scriptstyle [m,n]$}
\put(131,133){$\scriptstyle [m+1,n+1]$}

\end{picture}
\end{center}

\vspace{-1.3cm}

\noindent
where the vertices belong to the lattice $\ZZ^2$, the middle vertex represents
the weight $[m,n]$ of $w$ and the four vectors $(1,1)$, $(1,-1)$, $(-1,1)$
and $(-1,-1)$ indicate the action
of the respective lowering operators $z_{32}$, $z_{42}$, $z_{31}$ and $z_{41}$.
At most one step may be taken in any direction
beginning with $(1,1)$, then following with $(1,-1)$, $(-1,1)$
and $(-1,-1)$. The monomial \eqref{baskac} is obtained
by writing the product of the labells of the vectors used in the path and apply
the corresponding operator to $w$. For instance, the path of length two
consisting of $(1,-1)$ then following by $(-1,1)$ corresponds to the monomial
$z_{31}\tss z_{42}\tss w$ of weight $[m,n]$.

We thus obtain the following weight vectors in $K^+(m,n)$:
\begin{align}
[m,n] \qquad
& w,\quad z_{41}\tss z_{32}\tss w,\quad   z_{31}\tss z_{42}\tss w,\quad
z_{41}\tss z_{31}\tss z_{42}\tss z_{32}\tss w, \el\\
[m+1,n+1] \qquad
& z_{32}\tss w,\quad  z_{31}\tss z_{42}\tss z_{32}\tss w, \el\\
[m+1,n-1] \qquad
& z_{42}\tss w,\quad  z_{41}\tss z_{42}\tss z_{32}\tss w, \el\\
[m-1,n+1] \qquad
& z_{31}\tss w,\quad  z_{41}\tss z_{31}\tss z_{32}\tss w,\el\\
[m-1,n-1] \qquad
& z_{41}\tss w,\quad  z_{41}\tss z_{31}\tss z_{42}\tss w,\el\\
[m+2,n] \qquad
& z_{42}\tss z_{32}\tss w, \el\\
[m,n+2] \qquad
& z_{31}\tss z_{32}\tss w, \el\\
[m-2,n] \qquad
& z_{41}\tss z_{31}\tss w, \el\\
[m,n-2] \qquad
& z_{41}\tss z_{42}\tss w. \el
\end{align}

In the cases where $m\in\{0,1\}$ or $n\in\{0,1\}$
families of weight vectors in $K^+(m,n)$
are constructed
in a way similar to \eqref{baskac}.
We use the interpretation of the elements \eqref{baskac} as paths
in the same oriented graph with the additional condition that
all vertices $[k,l]$ of the path belong to the region
$k,l\geqslant 0$. We will call such paths {\em admissible}.
Clearly, if $m,n\geqslant 2$ then all paths are admissible.
For instance, in the case $m=n=0$ this leaves the elements
\begin{align}
[0,0] \qquad
& w,\quad
z_{41}\tss z_{31}\tss z_{42}\tss z_{32}\tss w, \el\\
[1,1] \qquad
& z_{32}\tss w,\quad  z_{31}\tss z_{42}\tss z_{32}\tss w, \el\\
[2,0] \qquad
& z_{42}\tss z_{32}\tss w, \el\\
[0,2] \qquad
& z_{31}\tss z_{32}\tss w, \el
\end{align}
of the space $K^+(0,0)$. In all remaining cases such vectors are listed
in Appendix~\ref{sec:abk}.

\bpr\label{prop:baskacgen}
A basis of the Kac module $K(m,n)$ is formed by
the vectors
\beql{baska}
E_{21}^k\tss E_{43}^l\ts z_{41}^{\theta_1}\ts
z_{31}^{\theta_2}\ts z_{42}^{\theta_3}\ts z_{32}^{\theta_4}\ts w,
\eeq
with the condition that the corresponding elements \eqref{baskac}
are associated with admissible paths, and where
\ben
k=0,1,\dots,m-\theta_1-\theta_2+\theta_3+\theta_4\Fand
l=0,1,\dots,n-\theta_1+\theta_2-\theta_3+\theta_4.
\een
\epr

\bpf
The restriction of the module $K(m,n)$ to the subalgebra $\kgot$
is a direct sum
\ben
K(m,n)\Big|_{\kgot}\cong \bigoplus_{r,s\geqslant 0} c_{r,s}\tss L^0(r,s),
\een
where the multiplicity $c_{r,s}$ is found by
\ben
c_{r,s}=\dim K^+(m,n)_{[r,s]},
\een
where the subscript $[r,s]$ indicates the corresponding weight subspace.
The subspace $K^+(m,n)_{[r,s]}$ coincides with the image of the
weight space $K(m,n)_{[r,s]}$ under the action of the
extremal projector $p=p(\kgot)$,
\ben
K^+(m,n)_{[r,s]}=p\tss K(m,n)_{[r,s]}.
\een
By the Poincar\'e--Birkhoff--Witt theorem,
the Kac module $K(m,n)$ is spanned by vectors of the form
\ben
E_{21}^k\tss E_{43}^l\tss E_{41}^{\theta_1}\tss
E_{31}^{\theta_2}\tss E_{42}^{\theta_3}\tss  E_{32}^{\theta_4}\tss w
\een
where each
$\theta_i$ takes values in $\{0,1\}$. Due to the properties \eqref{eppe}
of $p$, we may conclude that the space $K^+(m,n)_{[r,s]}$ is spanned
by the vectors
\beql{pevec}
p \tss E_{41}^{\theta_1}\tss
E_{31}^{\theta_2}\tss E_{42}^{\theta_3}\tss  E_{32}^{\theta_4}\tss w
\eeq
such that $r=m-\theta_1-\theta_2+\theta_3+\theta_4$ and
$s=n-\theta_1+\theta_2-\theta_3+\theta_4$. However, each vector \eqref{pevec}
is a linear combination of admissible elements of the form \eqref{baskac}.
Indeed, this follows by application of the formulas of Lemma~\ref{lem:inve}:
first replace $E_{32}$ with $z_{32}$, then use
the expression provided by Lemma~\ref{lem:inve} for $E_{42}$ to write the vector
as a linear combination of elements
$
p \tss E_{41}^{\theta_1}\tss
E_{31}^{\theta_2}\tss z_{42}^{\theta_3}\tss  z_{32}^{\theta_4}\tss w
$
and then use such replacements for $E_{31}$ and $E_{41}$.

Furthermore, each nonzero element of $K^+(m,n)_{[r,s]}$ generates
a $\kgot$-submodule of $K(m,n)$ of dimension $(r+1)(s+1)$.
Therefore, the module $K(m,n)$ is spanned by all
vectors \eqref{baska}. On the other hand, the number of these vectors
is easily calculated. For $m,n\geqslant 2$ it equals
\begin{align}
& 4\tss (m+1)(n+1)+ 2\tss (m+2)(n+2) +2\tss (m+2)n +2\tss m(n+2) +2\tss mn \el\\
&+(m+3)(n+1)+(m+1)(n+3)+(m-1)(n+1)+(m+1)(n-1)\el \\
&=16\tss (m+1)(n+1)\el
\end{align}
which coincides with $\dim K(m,n)$. This proves that the vectors
form a basis of $K(m,n)$. The same calculation in the cases where $m\leqslant 1$
or $n\leqslant 1$ confirms that the number of vectors matches $\dim K(m,n)$.
\epf

The proof of Proposition~\ref{prop:baskacgen} essentially contains
the decompositions of the Kac modules as $\kgot$-modules. In particular,
for $m,n\geqslant 2$ we have
\begin{multline}
K(m,n)\Big|_{\kgot}\cong 4\tss L^0(m,n)\oplus
2\tss L^0(m+1,n+1)\oplus 2\tss L^0(m+1,n-1)\\[0.3em]
\oplus 2\tss L^0(m-1,n+1)\oplus
2\tss L^0(m-1,n-1)\oplus L^0(m+2,n)\\[0.5em]
\oplus L^0(m,n+2)\oplus L^0(m-2,n)\oplus L^0(m,n-2).
\non
\end{multline}

Following the terminology used for representations
of simple Lie superalgebras \cite{k:ct}, we
will call the weight $[m,n]$ {\em typical}, if
the Kac module $K(m,n)$ is irreducible.
Otherwise, $[m,n]$ will be called {\em atypical}.
We will give necessary and sufficient conditions for
$[m,n]$ to be typical. They turn out to coincide with such conditions
for representations of $\psl(2|2)$ (see \cite{k:ct}, \cite{k:rc}), but the structure
of the atypical Kac modules differs; see also \cite{gqs:tp}.
Our main instrument will be the techniques of
Mickelsson--Zhelobenko algebras which will allow us to describe
$K(m,n)$ as a module over $\kgot$.

\bpr\label{prop:typ}
If $m\ne n$ then the Kac module $K(m,n)$ is irreducible.
\epr

\bpf
Observe that if $R$ is a nonzero submodule of $K(m,n)$, then the subspace
$R^+$ defined in \eqref{vplus} is a nonzero $S(\g,\kgot)$-submodule of $K^+(m,n)$.
Therefore, to describe $\g$-submodules of $K(m,n)$ it will be sufficient
to describe $S(\g,\kgot)$-submodules of $K^+(m,n)$.
Since $R^+$ is $\h$-invariant,
each weight component of $R^+$ is contained in $R^+$.
Working case by case for each weight subspace, we verify easily
with the use of formulas of Lemmas~\ref{lem:lowlow}, \ref{lem:rairai}, \ref{lem:railow}
and Appendix~\ref{sec:aar},
that the condition $R^+\ne\{0\}$
implies that $R^+$ contains the vector $w$.
For example, suppose that $m,n\geqslant 2$ and that a linear combination
\ben
c_1\tss w + c_2\tss  z_{41}\tss z_{32}\tss w +c_3\tss   z_{31}\tss z_{42}\tss w
+c_4\tss z_{41}\tss z_{31}\tss z_{42}\tss z_{32}\tss w,\qquad c_i\in\CC,
\een
belongs to $R^+$. Applying the operators $z_{23}$ and $z_{24}$
to this element, we obtain
the following two relations, respectively,
\ben
0= \frac{m-n}{2}\ts c_2 +\frac{m+n+2}{2\tss(n+1)}\ts c_3 \Fand
0= -\frac{m-n}{2\tss(n+2)}\ts c_2 +\frac{n\tss(m+n+2)}{2\tss(n+1)}\ts c_3
\een
together with $c_4=0$. The relations imply $c_2=c_3=0$ when $m\ne n$.
The same argument applied to the remaining weight subspaces implies that
$w\in R^+$ and so $R^+=K^+(m,n)$,
which proves that
$R=K(m,n)$. Thus, the corresponding Kac module $K(m,n)$ is irreducible.
\epf

Now suppose that $m=n$ and introduce the $\Sr(\g,\kgot)$-submodules
of $K^+(n,n)$ by
\ben
S^+_n=\Sr(\g,\kgot)\tss z_{32}\tss w\Fand T^+_n=\Sr(\g,\kgot)\tss z_{41}\tss w.
\een
The corresponding submodules $S_n$ and $T_n$ of $K(n,n)$ are then defined by
\beql{sntn}
S_n=\U(\kgot)\ts S^+_n\Fand T_n=\U(\kgot)\ts T^+_n.
\eeq

\bpr\label{prop:subst}
The $\g$-modules $S_n$ are irreducible for $n\geqslant 0$ and
the $\g$-modules $T_n$ are irreducible for $n\geqslant 1$. Moreover,
$T_0=\{0\}$ and
we have a $\g$-module isomorphism
\beql{stisom}
S_{n-1}\cong T_n,\qquad n\geqslant 1.
\eeq
\epr

\bpf
By using the formulas of Appendices~\ref{sec:amz}
and \ref{sec:aar}, we can produce explicit bases of
$S^+_n$ and $T^+_n$. Arranging the basis vectors in accordance with their $\h$-weights,
for $S^+_n$ with $n\geqslant 1$ we have
\begin{align}
[n,n] \qquad
& z_{41}\tss z_{32}\tss w,\quad
z_{41}\tss z_{31}\tss z_{42}\tss z_{32}\tss w, \el\\
[n+1,n+1] \qquad
& z_{32}\tss w,\quad  z_{31}\tss z_{42}\tss z_{32}\tss w, \el\\
[n+1,n-1] \qquad
& z_{41}\tss z_{42}\tss z_{32}\tss w, \el\\
[n-1,n+1] \qquad
& z_{41}\tss z_{31}\tss z_{32}\tss w,\el\\
[n+2,n] \qquad
& z_{42}\tss z_{32}\tss w, \el\\
[n,n+2] \qquad
& z_{31}\tss z_{32}\tss w, \el
\end{align}
and for $T^+_n$ with $n\geqslant 2$ we have
\begin{align}
[n,n] \qquad
& z_{32}\tss z_{41}\tss w,\quad
z_{31}\tss z_{42}\tss z_{32}\tss z_{41}\tss w, \el\\
[n+1,n-1] \qquad
& z_{42}\tss z_{32}\tss z_{41}\tss w, \el\\
[n-1,n+1] \qquad
& z_{31}\tss z_{32}\tss z_{41}\tss w, \el\\
[n-1,n-1] \qquad
& z_{41}\tss w,\quad  z_{31}\tss z_{42}\tss z_{41}\tss w,\el\\
[n-2,n] \qquad
& z_{31}\tss z_{41}\tss w, \el\\
[n,n-2] \qquad
& z_{42}\tss z_{41}\tss w. \el
\end{align}
Similarly, the basis of $S^+_0$ is given by
\begin{align}
[0,0] \qquad
& w,\quad
z_{41}\tss z_{31}\tss z_{42}\tss z_{32}\tss w, \el\\
[1,1] \qquad
& z_{32}\tss w,\quad  z_{31}\tss z_{42}\tss z_{32}\tss w, \el\\
[2,0] \qquad
& z_{42}\tss z_{32}\tss w, \el\\
[0,2] \qquad
& z_{31}\tss z_{32}\tss w, \el
\end{align}
while the basis of $T^+_1$ is
\begin{align}
[1,1] \qquad
& z_{32}\tss z_{41}\tss w,\quad
z_{31}\tss z_{42}\tss z_{32}\tss z_{41}\tss w, \el\\
[2,0] \qquad
& z_{42}\tss z_{32}\tss z_{41}\tss w, \el\\
[0,2] \qquad
& z_{31}\tss z_{32}\tss z_{41}\tss w, \el\\
[0,0] \qquad
& z_{41}\tss w,\quad  z_{41}\tss z_{31}\tss z_{42}\tss w,\el
\end{align}
and $T^+_0=\{0\}$.
Using formulas of Appendices~\ref{sec:amz}
and \ref{sec:aar} once again, we can see that the nonzero submodules
$S^+_n$ and $T^+_n$ of the $\Sr(\g,\kgot)$-module
$K^+(n,n)$ are irreducible. This implies that the corresponding submodules
$S_n$ and $T_n$ of $K(n,n)$ are also irreducible.

Finally, to prove the last statement of the proposition, for a given $n\geqslant 1$
denote by $w'$ the highest vector of the $\kgot$-module $L^0(n-1,n-1)$.
The Mickelsson--Zhelobenko algebra relations
imply that for $n\geqslant 2$ we have an $\Sr(\g,\kgot)$-module isomorphism
\ben
\phi:S^+_{n-1}\to T^+_n,\qquad z_{32}\tss w'\mapsto z_{32}\tss z_{41}\tss w,
\een
with the inverse map given by
\ben
\phi^{-1}: z_{41}\tss w\mapsto -\frac{n+1}{n}\ts z_{41}\tss z_{32}\tss w'.
\een
For $n=1$ the statement is equivalent to the existence
of an isomorphism $K^+(0,0)\cong T^+_1$. It is provided by the map
\ben
\phi:K^+(0,0)\to T^+_1,\qquad w'\mapsto z_{41}\tss w.
\een
This yields the desired isomorphism \eqref{stisom}.
\epf

\bpr\label{prop:atyp}
The Kac module $K(n,n)$ over $\g$ with $n\geqslant 1$ is the direct sum of two
irreducible submodules,
\ben
K(n,n)=S_n\oplus T_n.
\een
The module $K(0,0)=S_0$ is irreducible. Hence, we have an isomorphism
\ben
K(n,n)\cong S_n\oplus S_{n-1},\qquad n\geqslant 0,
\een
assuming $S_{-1}=\{0\}$.
\epr

\bpf
This will follow from Proposition~\ref{prop:subst}. It suffices to verify that
\beql{kpst}
K^+(n,n)=S^+_n\oplus T^+_n,\qquad n\geqslant 1.
\eeq
However, $\dim K^+(n,n)=\dim S^+_n+\dim T^+_n$ and we have
\ben
K^+(n,n)=S^+_n+ T^+_n
\een
due to the relation
\ben
w=
-\frac{n+2}{n+1}\ts z_{41}\tss z_{32} w
-\frac{n+1}{n}\ts z_{32}\tss z_{41}\tss w;
\een
see Lemma~\ref{lem:lowlow}.
Therefore, the intersection of $S^+_n$ and $T^+_n$ is zero
and \eqref{kpst} follows, thus completing the proof.
\epf

\subsection{Classification theorem}
\label{subsec:ct}

We can now prove the classification theorem for representations
of the Lie superalgebra $\g=\psl(2|2)\ltimes\CC^3$, where the central elements act
by
\beql{ckp}
C\mapsto 0,\qquad K\mapsto 0\Fand P\mapsto 1.
\eeq

\bth\label{thm:class}
A complete list of pairwise non-isomorphic finite-dimensional
irreducible representations of $\g$ with the conditions \eqref{ckp} consists of
\begin{enumerate}
\item
the Kac modules $K(m,n)$ with $m,n\in\ZZ_+$ and $m\ne n$,\newline
$\dim K(m,n)=16\tss(m+1)(n+1)$,
\item
the modules $S_n$ with $n\in\ZZ_+$,\quad $\dim S_n=8\tss(n+1)(n+2)$.
\end{enumerate}
\eth

\bpf
Consider the following triangular decomposition of the Lie superalgebra $\g$,
\ben
\g=\bar\n^-\oplus\h\oplus\bar\n^+,
\een
where $\h$ is spanned by the elements $h_1$, $h_3$, $C$, $K$ and $P$,
whereas the subalgebras $\bar\n^+$
and $\bar\n^-$ are defined by
\ben
\bal
\bar\n^+&=\text{span of }\{E_{12}, E_{34}, E_{31}, E_{32}, E_{14}, E_{24}\},\\
\bar\n^-&=\text{span of }\{E_{21}, E_{43}, E_{13}, E_{23}, E_{41}, E_{42}\}.
\eal
\een
Given a pair of complex numbers $\mu=(\mu_1,\mu_3)$, consider the one-dimensional
representation $\CC_{\mu}$ of the Lie superalgebra $\h\oplus \bar\n^+$
defined by
\ben
\bar\n^+\ts 1_{\mu}=0,\quad h_1\ts 1_{\mu}=\mu_1 \ts 1_{\mu},\quad
h_3\ts 1_{\mu}=\mu_3 \ts 1_{\mu},\quad C\ts 1_{\mu}=0,
\quad K\ts 1_{\mu}=0,\quad P\ts 1_{\mu}=1_{\mu},
\een
where $1_{\mu}$ denotes the basis vector of $\CC_{\mu}$.
The corresponding {\em Verma module} $\overline M(\mu)$ is then defined by
\ben
\overline M(\mu)=\U(\g)\ot_{U(\h\oplus \bar\n^+)}\CC_{\mu}.
\een
By a standard argument, $\overline M(\mu)$ contains a unique maximal
proper submodule $V$ and we set $\overline L(\mu)=\overline M(\mu)/V$.
The module $\overline L(\mu)$ is irreducible and any finite-dimensional
irreducible representation of $\g$ is isomorphic to $\overline L(\mu)$
for a certain uniquely determined $\mu$. Therefore, to classify
finite-dimensional irreducible representations of $\g$ it suffices to find
necessary and sufficient conditions on $\mu$ for the module $\overline L(\mu)$
to be finite-dimensional.

By considering the $\U(\kgot)$-cyclic span of the vector $1_{\mu}$ we conclude
that the condition $\dim \overline L(\mu)<\infty$ implies that both
$\mu_1$ and $\mu_3$ are nonnegative integers. In what follows we will
assume that $\mu_1,\mu_3\in\ZZ_+$. As a next step, we will demonstrate
that $\overline L(\mu)$ is infinite-dimensional unless $\mu_3\geqslant 2$.

Suppose first that $\mu_3=0$. The vector
\ben
v=z_{42}\tss 1_{\mu}=E_{42}\tss 1_{\mu}\in \overline L(\mu)
\een
is nonzero since $E_{31}\ts v=P\ts 1_{\mu}=1_{\mu}$. On the other hand,
the weight of $v$ is $(\mu_1+1,\ts -1)$ and $\kgot^+\ts v=0$.
Therefore, the $\kgot$-module $\U(\kgot)\ts v$ is infinite-dimensional
and so is $\overline L(\mu)$.

Now let $\mu_3=1$. If $\mu_1\geqslant 1$ then the vector
\ben
u=z_{41}\tss z_{42}\tss 1_{\mu}=E_{41}\tss E_{42}\tss 1_{\mu}-P\tss
E_{43}\tss 1_{\mu}=E_{41}\tss E_{42}\tss 1_{\mu}-
E_{43}\tss 1_{\mu}\in \overline L(\mu)
\een
is nonzero since
\ben
E_{24}\tss E_{14}\tss u=-\frac{\mu_1(\mu_1+1)}{2}\ts 1_{\mu}
\een
while the weight of $u$ is $(\mu_1,-1)$. If $\mu_1=0$ then the vector
$y=E_{41}\tss 1_{\mu}$ is nonzero since $E_{32}\tss y=-P\tss 1_{\mu}=-1_{\mu}$
and the weight of $y$ is $(-1,0)$.

As a final step, we will show that each representation $\overline L(\mu)$
with $\mu_3\geqslant 2$
is finite-dimensional and isomorphic to exactly one module in the list
given in the formulation of the theorem.
By the construction of the
Kac modules $K(m,n)$, for $\mu=(m,n+2)$ we have the homomorphism
\ben
\overline M(m,n+2)\to K(m,n),\qquad 1_{\mu}\mapsto z_{31}\tss z_{32}\tss w.
\een
By Propositions~\ref{prop:typ} and \ref{prop:subst}, this yields
isomorphisms
\ben
\overline L(m,n+2)\cong K(m,n),\qquad m,n\in\ZZ_+,\quad m\ne n,
\een
and
\ben
\overline L(n,n+2)\cong S_n,\qquad n\in\ZZ_+.
\een
The theorem is proved.
\epf

The techniques of Mickelsson--Zhelobenko algebras which we used to prove
Theorem~\ref{thm:class} can also be applied to reproduce the
well-known descriptions of finite-dimensional
irreducible representations of the Lie superalgebras $\psl(2|2)$ and $\sll(2|2)$
\cite{k:rc}; see also \cite{kkp:fd}, \cite{ps:fd} and \cite{zg:uc}.
Namely, the above arguments can be easily modified
for the case where the central element $P$ of $\g$ acts as the zero operator.
The corresponding Kac modules $K^{\circ}(m,n)$ over $\psl(2|2)$ are irreducible
for $m\ne n$. However, the structure of $K^{\circ}(n,n)$ differs from that of
the $\g$-module $K(n,n)$; the corresponding submodules $S^{\circ}_n$
and $T^{\circ}_n$ of $K^{\circ}(n,n)$, which are defined exactly as in \eqref{sntn},
are not irreducible for $n\geqslant 1$.
The intersection
$U_n=S^{\circ}_n\cap T^{\circ}_n$ is nonzero for $n\geqslant 1$ and we have
\ben
U^+_n=\text{span of }\{z_{41}\tss z_{32}\tss w,\quad
z_{41}\tss z_{31}\tss z_{42}\tss z_{32}\tss w,\quad
z_{41}\tss z_{42}\tss z_{32}\tss w,\quad
z_{41}\tss z_{31}\tss z_{32}\tss w\}.
\een
The sum $R_n=S^{\circ}_n+T^{\circ}_n$
is a proper submodule of $K^{\circ}(n,n)$ (cf. Proposition~\ref{prop:subst})
and the quotient $L^{\circ}(n,n)=K^{\circ}(n,n)/R_n$ is irreducible.
The vectors
$
w,\ts
z_{31}\tss w,\ts
z_{42}\tss w,\ts
z_{31}\tss z_{42}\tss w
$
form a basis of $L^{\circ+}(n,n)$ for $n\geqslant 1$.

To summarize,
we get the following description of $\psl(2|2)$-modules;
see \cite{k:ct} and \cite{k:rc}.

\bco\label{cor:pze}
A complete list of pairwise non-isomorphic finite-dimensional
irreducible representations of $\psl(2|2)$ consists of
\begin{enumerate}
\item
the Kac modules $K^{\circ}(m,n)$ with $m,n\in\ZZ_+$ and
$m\ne n$,\newline $\dim K^{\circ}(m,n)=16\tss(m+1)(n+1)$,
\item
the modules $L^{\circ}(n,n)$ with $n\geqslant 1$,\quad $\dim L^{\circ}(n,n)
=4\tss n(n+2)+2$,
\item
the trivial one-dimensional module $L^{\circ}(0,0)$.
\qed
\end{enumerate}
\eco

To state the corresponding results for the Lie superalgebra
$\g'=\sll(2|2)=\psl(2|2)\oplus \CC C$, consider the Kac modules $K(m,n;2\tss c)$
over $\g'$ (with $m,n\in\ZZ_+$),
which are defined as in \eqref{kacmodule}, except that
the central elements now act by
$C\mapsto c$, $K\mapsto 0$ and $P\mapsto 0$, and we assume that
the complex number $c$ is nonzero.
The $\g'$-module $K(m,n;2\tss c)$ is irreducible if and only if
\beql{typco}
m-n\ne \pm \tss 2\tss c\Fand m+n+2\ne \pm\tss 2\tss c.
\eeq
These conditions define the class of {\em typical} representations of $\g'$.
The remaining {\em atypical} representations are nontrivial
quotients of $K(m,n;2\tss c)$ in the cases where \eqref{typco} does not hold.
To describe the corresponding submodules,
consider the Mickelsson algebra $\Sr(\g',\kgot)$
and introduce submodules
\begin{alignat}{2}
S&\subset K(m,n;m-n),\qquad &&T\subset K(m,n;n-m),
\non\\
X&\subset K(m,n;-m-n-2),\qquad &&Y\subset  K(m,n;m+n+2),
\non
\end{alignat}
by setting $V=\U(\kgot)\tss V^+$, where $V$ denotes one of the four submodules,
and $V^+$ is the $\Sr(\g',\kgot)$-submodule of the respective Kac module,
\begin{alignat}{2}
S^+&=\Sr(\g',\kgot)\tss z_{32}\tss w,\qquad &&T^+=\Sr(\g',\kgot)\tss z_{41}\tss w
\non\\
X^+&=\Sr(\g',\kgot)\tss z_{31}\tss w,\qquad &&Y^+=\Sr(\g',\kgot)\tss z_{42}\tss w.
\non
\end{alignat}

\bco\label{cor:sl22}
A complete list of pairwise non-isomorphic finite-dimensional
irreducible representations of $\alg{sl}(2|2)$ with a non-zero
eigenvalue $c$ of the central element $C$
consists of
\begin{enumerate}
\item the Kac modules $K(m,n;2\tss c)$ with the conditions \eqref{typco}, \\
$\dim K(m,n;2\tss c)=16(m+1)(n+1)$,
\item the modules $K(m,n;m-n)/S$, \\
$\dim K(m,n;m-n)/S=4\big(m(n+1)+(m+1)n\big)$,
\item the modules $K(m,n;-m+n)/T$, \\
$\dim K(m,n;-m+n)/T=4\big((m+1)(n+2)+(m+2)(n+1)\big)$,
\item the modules $K(m,n;-m-n-2)/X$, \\
$\dim K(m,n;-m-n-2)/X=4\big((m+2)(n+1)+(m+1)n\big)$,
\item the modules $K(m,n;m+n+2)/Y$, \\
$\dim K(m,n;-m-n-2)/Y=4\big((m+1)(n+2)+m(n+1)\big)$.
\qed
\end{enumerate}
\eco

\subsection{Explicit construction of representations}
\label{subsec:ec}

Our proof of the classification theorem (Theorem~\ref{thm:class}) was based on explicit bases
of irreducible representations $V$. They all have the form
\beql{baskaefo}
E_{21}^k\tss E_{43}^l\ts z_{41}^{\theta_1}\ts
z_{31}^{\theta_2}\ts z_{42}^{\theta_3}\ts z_{32}^{\theta_4}\ts w,
\eeq
with some conditions on the parameters, where $w$ is the highest vector
of the $\kgot$-module $L^0(m,n)$.
The matrix elements
for the action of the generators of
$\g$ in this basis can be found from the Mickelsson--Zhelobenko
algebra relations in a standard way; cf.\cite[Ch.~9]{m:yc} and \cite{m:cb}.
First observe that if
$k>m-\theta_1-\theta_2+\theta_3+\theta_4$ or
$l>n-\theta_1+\theta_2-\theta_3+\theta_4$ then the corresponding vector
\eqref{baskaefo} is zero. This is easily verified by considering all
possible values of the parameters $\theta_i$. For example, if
$\theta_1=\theta_3=\theta_4=0$ and $\theta_2=1$ then for $m\geqslant 1$ we have
\ben
\bal
E_{21}^m\tss z_{31}\tss w
&=E_{21}^m\Big(E_{31}+E_{21}\tss E_{32}\ts\frac1{m+1}\Big)\tss w\\
{}&=E_{31}\tss E_{21}^m\tss w+\frac1{m+1}\ts E_{32}\tss E_{21}^{m+1}\tss w
-(m+1)\tss E_{31}\tss E_{21}^m\ts \frac1{m+1}\tss w=0
\eal
\een
since $E_{21}^{m+1}\tss w=0$ in $K(m,n)$. Therefore, $E_{21}^k\tss z_{31}\tss w=0$
for all $k>m$. This determines the action of the subalgebra $\kgot$
on the basis vectors.

Furthermore,
using the commutation relations of $\g$ we can reduce
the calculation to the case where one of the generators of the form
$E_{ik}$ or $E_{ki}$ with $i=1,2$ and $k=3,4$ acts on the vector $v\in V^+$.
Then we write this generator in the form provided by
Lemma~\ref{lem:inve} and apply the formulas for the action
of the elements $z_{ik}$ and $z_{ki}$ on the basis of $V^+$.
To illustrate, consider the action of $E_{13}$
on the basis vector \eqref{baskaefo},
\ben
E_{13}\ts
E_{21}^k\tss E_{43}^l\ts z_{41}^{\theta_1}\ts
z_{31}^{\theta_2}\ts z_{42}^{\theta_3}\ts z_{32}^{\theta_4}\ts w
=E_{21}^k\tss E_{43}^l\tss E_{13}\ts z_{41}^{\theta_1}\ts
z_{31}^{\theta_2}\ts z_{42}^{\theta_3}\ts z_{32}^{\theta_4}\ts w
-k\ts E_{21}^{k-1} E_{43}^l\tss E_{23}\ts z_{41}^{\theta_1}\ts
z_{31}^{\theta_2}\ts z_{42}^{\theta_3}\ts z_{32}^{\theta_4}\ts w.
\een
Next, replace $E_{13}$ and $E_{23}$ by their expression
provided by Lemma~\ref{lem:inve}. In particular, the first vector on the right hand side
becomes
\ben
E_{21}^k\tss E_{43}^l\tss
\Big(z_{13}-E_{43}\tss z_{14}\frac{1}{h_3+1}\Big)\ts z_{41}^{\theta_1}\ts
z_{31}^{\theta_2}\ts z_{42}^{\theta_3}\ts z_{32}^{\theta_4}\ts w
\een
so that the calculation is completed by applying the formulas of Appendix~\ref{sec:aar}
for the action of $z_{13}$ and $z_{14}$ on the vector
$z_{41}^{\theta_1}\ts
z_{31}^{\theta_2}\ts z_{42}^{\theta_3}\ts z_{32}^{\theta_4}\ts w$.

We will omit explicit matrix element formulas
to avoid significant extension of the paper for the reason
that their reproduction is straightforward from the formulas
of Appendices~\ref{sec:amz} and \ref{sec:aar}.

Note that the same techniques of Mickelsson--Zhelobenko
algebras can also be used to reproduce explicit basis
constructions for representations of the Lie superalgebras $\sll(2|2)$ and $\gl(2|2)$
given in \cite{kkp:fd}, \cite{ps:fd} and \cite{zg:uc}.

\appendix

\section{Relations in the Mickelsson--Zhelobenko algebra}
\label{sec:amz}
\setcounter{equation}{0}

The following relations in $\Z(\g,\kgot)$ are written without a specialization
of the values of the central elements $C$, $K$ and $P$.

\ble\label{lem:lowlow}
We have the relations for the lowering operators:
\begin{align}
z_{31}z_{41}&=-z_{41}z_{31}\frac{h_3}{h_3+1}, \non\\
z_{42}z_{41}&=-z_{41}z_{42}\frac{h_1}{h_1+1},  \non\\
z_{32}z_{41}
&=-P\frac{h_3}{h_3+1}- z_{31}z_{42}\frac{h_1-h_3}{(h_1+1)(h_3+1)}
- z_{41}z_{32}\frac{h_3(h_3+2)}{(h_3+1)^2},
\non\\
z_{42}z_{31}&=P-z_{31}z_{42} +z_{41}z_{32}\frac{h_1+h_3+2}{(h_1+1)(h_3+1)},\non\\
z_{32}z_{31}&= -z_{31}z_{32}\frac{h_1}{h_1+1}, \non\\
z_{32}z_{42}&= -z_{42}z_{32}\frac{h_3}{h_3+1}.\non
\end{align}
\ele

\ble\label{lem:rairai}
We have the relations for the raising operators:
\begin{align}
z_{14}z_{13}&=-z_{13}z_{14}\frac{h_3}{h_3+1}, \non\\
z_{14}z_{24}&=-z_{24}z_{14}\frac{h_1}{h_1+1},  \non\\
z_{14}z_{23}&= -K\frac{h_3}{h_3+1}-z_{24}z_{13}
\frac{h_1-h_3}{(h_1+1)(h_3+1)}-z_{23}z_{14}\frac{h_3(h_3+2)}{(h_3+1)^2},\non\\
z_{13}z_{24}&=K-z_{24}z_{13} +z_{23}z_{14}\frac{h_1+h_3+2}{(h_1+1)(h_3+1)},\non\\
z_{13}z_{23}&= -z_{23}z_{13}\frac{h_1}{h_1+1}, \non\\
z_{24}z_{23}&= -z_{23}z_{24}\frac{h_3}{h_3+1}.\non
\end{align}
\ele

\ble\label{lem:railow}
We have the relations for the raising and lowering operators;

with $z_{14}$:
\begin{align}
z_{14}z_{41}&= \frac{h_1h_3(h_1-h_3+2\tss C)}{2(h_1+1)(h_3+1)} +z_{31}z_{13}
\frac{h_1(h_1+2)}{(h_1+1)^2(h_3+1)}-z_{32}z_{23}\frac{1}{(h_1+1)(h_3+1)} \el \\
&\quad -z_{41}z_{14}
\frac{h_1h_3(h_1+2)(h_3+2)}{(h_1+1)^2(h_3+1)^2}
+ z_{42}z_{24}\frac{h_3(h_3+2)}{(h_1+1)(h_3+1)^2}, \non\\[0.5em]
z_{14}z_{31}&=-z_{31}z_{14}\frac{h_1(h_1+2)}
{(h_1+1)^2}+z_{32}z_{24}\frac{1}{h_1+1}, \non\\[0.5em]
z_{14}z_{42}&=-z_{42}z_{14}\frac{h_3(h_3+2)}
{(h_3+1)^2}+z_{32}z_{13}\frac{1}{h_3+1}, \non\\[0.5em]
z_{14}z_{32}&=-z_{32}z_{14}, \non
\end{align}

with $z_{13}$:
\begin{align}
z_{13}z_{41}&=-z_{41}z_{13}\frac{h_1(h_1+2)}{(h_1+1)^2}
+z_{42}z_{23}\frac{1}{h_1+1}, \non\\[0.5em]
z_{13}z_{31}&= \frac{h_1(h_1+h_3+2+2\tss C)}{2(h_1+1)}
-z_{31}z_{13}\frac{h_1(h_1+2)}{(h_1+1)^2} +z_{32}z_{23}\frac{1}{h_1+1} \non\\
&\quad - z_{41}z_{14}\frac{h_1(h_1+2)}{(h_1+1)^2(h_3+1)}
+z_{42}z_{24}\frac{1}{(h_1+1)(h_3+1)},\non\\[0.5em]
z_{13}z_{42}&=-z_{42}z_{13}, \non\\[0.5em]
z_{13}z_{32}&=-z_{32}z_{13}-z_{42}z_{14}\frac{1}{h_3+1}, \non
\end{align}

with $z_{24}$:
\begin{align}
z_{24}z_{41}&=-z_{41}z_{24}\frac{h_3(h_3+2)}{(h_3+1)^2}
+ z_{31}z_{23}\frac{1}{h_3+1},  \non\\[0.5em]
z_{24}z_{31}&=-z_{31}z_{24}, \non\\[0.5em]
z_{24}z_{42}&= -\frac{h_3(h_1+h_3+2-2\tss C)}{2(h_3+1)} +z_{31}z_{13}
\frac{1}{(h_1+1)(h_3+1)} +z_{32}z_{23}\frac{1}{h_3+1} \non\\
&\quad - z_{41}z_{14}\frac{h_3(h_3+2)}{(h_1+1)(h_3+1)^2}
-z_{42}z_{24}\frac{h_3(h_3+2)}{(h_3+1)^2}, \non\\[0.5em]
z_{24}z_{32}&=-z_{32}z_{24}-z_{31}z_{14}\frac{1}{h_1+1}, \non
\end{align}

with $z_{23}$:
\begin{align}
z_{23}z_{41}&=-z_{41}z_{23}, \non\\
z_{23}z_{31}&=-z_{31}z_{23}-z_{41}z_{24}\frac{1}{h_3+1},  \non\\
z_{23}z_{42}&=-z_{42}z_{23}-z_{41}z_{13}\frac{1}{h_1+1},  \non\\
z_{23}z_{32}&= -\frac{h_1-h_3-2\tss C}{2} -z_{31}z_{13}\frac{1}{h_1+1}-z_{32}z_{23}
-z_{41}z_{14}\frac{1}{(h_1+1)(h_3+1)} - z_{42}z_{24}\frac{1}{h_3+1}.\non
\end{align}
\ele

\section{Action of raising operators}
\label{sec:aar}
\setcounter{equation}{0}

Lemmas~\ref{lem:lowlow}, \ref{lem:rairai} and \ref{lem:railow} imply the following
relations for the action of the raising operators on the vectors
of the space $K^+(m,n)$ of $\kgot^+$-invariants of the Kac module $K(m,n)$.

\bigskip
\noindent
Action on $w$:

\ben
z_{14}\cdot w=z_{24}\cdot w=z_{13}\cdot w=z_{23}\cdot w=0.
\een

\noindent
Action on $z_{41}z_{32}\tss w$:
\begin{align}
&z_{14}\cdot z_{41}z_{32}\tss w=\frac{(m-n)
[(m+1)(n+1)+1]}{2(m+2)(n+2)} z_{32}\tss w,\el\\[0.3em]
&z_{24}\cdot z_{41}z_{32}\tss w=-\frac{m-n}{2(n+2)}z_{31}\tss w, \el \\[0.3em]
&z_{13}\cdot z_{41}z_{32}\tss w=-\frac{m-n}{2(m+2)}z_{42}\tss w, \el\\[0.3em]
&z_{23}\cdot z_{41}z_{32}\tss w=\frac{m-n}{2}z_{41}\tss w. \el
\end{align}

\noindent
Action on $z_{31}z_{42}\tss w$:
\begin{align}
&z_{14}\cdot z_{31}z_{42}\tss w=-\frac{n(m+n+2)}
{2(m+2)(n+1)}z_{32}\tss w, \el\\[0.3em]
&z_{24}\cdot z_{31}z_{42}\tss w= \frac{n(m+n+2)}
{2(n+1)}z_{31}\tss w, \el\\[0.3em]
&z_{13}\cdot z_{31}z_{42}\tss w= \frac{(mn+m+n)
(m+n+2)}{2(m+2)(n+1)}z_{42}\tss w, \el\\[0.3em]
&z_{23}\cdot z_{31}z_{42}\tss w= \frac{m+n+2}{2(n+1)}z_{41}\tss w.
\el
\end{align}

\noindent
Action on $z_{41}z_{31}z_{42}z_{32}\tss w$:
\begin{align}
&z_{14}\cdot z_{41}z_{31}z_{42}z_{32}\tss w= \frac{m+n+2}{2(m+2)}Pz_{32}\tss w
+\frac{(m-n)(m+1)(mn+m+2n+1)}{2(m+2)^2(n+2)}\el \\
&\hspace{5cm} \times \Bigl[ 1
+\frac{(m+n+2)(mn+m+2n+3)}{(m+1)
(n+1)(mn+m+2n+1)} \Bigr] z_{31}z_{42}z_{32}\tss w, \el \\[0.5em]
&z_{24}\cdot z_{41}z_{31}z_{42}z_{32}\tss w= -\frac{(m+n+2)
(n+1)(n^2+5n+7)}{2(n+2)^3} z_{41}z_{31}z_{32}\tss w, \el\\[0.5em]
&z_{13}\cdot z_{41}z_{31}z_{42}z_{32}\tss w= -\frac{(m+n+2)
(mn+2m+n+3)}{2(m+2)^2} z_{41}z_{42}z_{32}\tss w
-\frac{m-n}{2(m+2)}Pz_{42}\tss w, \el\\[0.5em]
&z_{23}\cdot z_{41}z_{31}z_{42}z_{32}\tss w= \frac{m-n}{2}z_{41}z_{31}z_{42}\tss w.
\el
\end{align}

\noindent
Action on $z_{32}\tss w$:
\begin{align}
& z_{23}\cdot z_{32}\tss w=-\frac{m-n}{2} \tss w, \qquad
z_{14}\cdot z_{32}\tss w=z_{24}\cdot z_{32}\tss w=z_{13}\cdot z_{32}\tss w=0.\el
\end{align}

\noindent
Action on $z_{31}z_{42}z_{32}\tss w$:
\begin{align}
&z_{14}\cdot z_{31}z_{42}z_{32}\tss w=0, \el \\
&z_{24}\cdot z_{31}z_{42}z_{32}\tss w=\frac{m+n+2}{2} z_{31}z_{32}\tss w, \el\\
&z_{13}\cdot z_{31}z_{42}z_{32}\tss w=\frac{m+n+2}{2} z_{42}z_{32}\tss w,   \el \\
&z_{23}\cdot z_{31}z_{42}z_{32}\tss w=
-\frac{m-n}{2} z_{31}z_{42}\tss w+
\frac{m+n+2}{2(n+1)} z_{41}z_{32}\tss w.
\el
\end{align}

\noindent
Action on $z_{42}\tss w$:
\begin{align}
&z_{24}\cdot z_{42}\tss w=-\frac{n(m+n+2)}{2(n+1)}\tss w, \qquad
z_{14}\cdot z_{42}\tss w=z_{13}\cdot z_{42}\tss w=z_{23}\cdot z_{42}\tss w=0.\el
\end{align}

\noindent
Action on $z_{41}z_{42}z_{32}\tss w$:
\begin{align}
&z_{14}\cdot z_{41}z_{42}z_{32}\tss w=
\frac{n(m-n)}{2(n+1)} z_{42}z_{32}\tss w,\el\\
&z_{24}\cdot z_{41}z_{42}z_{32}\tss w=
\frac{n(n+2)(m+n+2)}{2(n+1)^2}z_{41}z_{32}\tss w
+\frac{m-n}{2(n+1)}z_{31}z_{42}\tss w,  \el\\
&z_{13}\cdot z_{41}z_{42}z_{32}\tss w= 0,\el \\
&z_{23}\cdot z_{41}z_{42}z_{32}\tss w= -\frac{m-n}{2}z_{41}z_{42}\tss w.
\el
\end{align}

\noindent
Action on $z_{31}\tss w$:
\begin{align}
&z_{13}\cdot z_{31}\tss w=\frac{m(m+n+2)}{2(m+1)}\tss w, \qquad
z_{14}\cdot z_{31}\tss w=z_{24}\cdot z_{31}\tss w=z_{23}\cdot z_{31}\tss w=0. \el
\end{align}

\noindent
Action on $z_{41}z_{31}z_{32}\tss w$:
\begin{align}
&z_{14}\cdot z_{41}z_{31}z_{32}\tss w=
\frac{m(m-n)}{2(m+1)} z_{31}z_{32}\tss w,\el\\
&z_{24}\cdot z_{41}z_{31}z_{32}\tss w= 0,\el\\
&z_{13}\cdot z_{41}z_{31}z_{32}\tss w=
-\frac{m(m+n+2)}{2(m+1)} z_{41}z_{32}\tss w
-\frac{m-n}{2(m+1)} (z_{31}z_{42}\tss w-P\tss w),\el\\
&z_{23}\cdot z_{41}z_{31}z_{32}\tss w= -\frac{m-n}{2}z_{41}z_{31}\tss w.
\el
\end{align}

\noindent
Action on $z_{41}\tss w$:
\begin{align}
&z_{14}\cdot z_{41}\tss w=
\frac{mn(m-n)}{2(m+1)(n+1)}\tss w,  \qquad
z_{24}\cdot z_{41}\tss w=z_{13}\cdot z_{41}\tss w=z_{23}\cdot z_{41}\tss w=0.  \el
\end{align}

\noindent
Action on $z_{41}z_{31}z_{42}\tss w$:
\begin{align}
&z_{14}\cdot z_{41}z_{31}z_{42}\tss w=
\frac{(m-n)(mn+m+n+2)}{2(m+1)(n+1)}z_{31}z_{42}\tss w
+\frac{(m+n+2)n(n+2)}
{2(m+1)(n+1)^2}z_{41}z_{32}\tss w\el\\
&\qquad\qquad\qquad\qquad +\frac{(m+n+2)n}{2(m+1)}P \tss w, \el\\
&z_{24}\cdot z_{41}z_{31}z_{42}\tss w=
-\frac{n(m+n+2)}{2(n+1)}z_{41}z_{31}\tss w, \el\\
&z_{13}\cdot z_{41}z_{31}z_{42}\tss w=
-\frac{m(m+n+2)}{2(m+1)}z_{41}z_{42}\tss w, \el\\
&z_{23}\cdot z_{41}z_{31}z_{42}\tss w= 0.
\el
\end{align}

\noindent
Action on $z_{42}z_{32}\tss w$:
\begin{align}
&z_{14}\cdot z_{42}z_{32}\tss w=0, && z_{24}\cdot
z_{42}z_{32}\tss w=\frac{m+n+2}{2}z_{32}\tss w, \el\\
&z_{13}\cdot z_{42}z_{32}\tss w=0, && z_{23}\cdot
z_{42}z_{32}\tss w=\frac{m-n}{2}z_{42}\tss w.
\el
\end{align}

\noindent
Action on $z_{31}z_{32}\tss w$:
\begin{align}
&z_{14}\cdot z_{31}z_{32}\tss w=0, && z_{24}\cdot z_{31}z_{32}\tss w=0, \el\\
&z_{13}\cdot z_{31}z_{32}\tss w=
\frac{m+n+2}{2}z_{32}\tss w, && z_{23}
\cdot z_{31}z_{32}\tss w=\frac{m-n}{2}z_{31}\tss w.
\el
\end{align}

\noindent
Action on $z_{41}z_{31}\tss w$:
\begin{align}
&z_{14}\cdot z_{41}z_{31}\tss w=\frac{(m-1)
(m-n)}{2m}z_{31}\tss w, && z_{24}\cdot z_{41}z_{31}\tss w=0, \el\\
&z_{13}\cdot z_{41}z_{31}\tss w=-\frac{(m-1)(m+n+2)}{2m}
z_{41}\tss w, && z_{23}\cdot z_{41}z_{31}\tss w=0.
\el
\end{align}

\noindent
Action on $z_{41}z_{42}\tss w$:
\begin{align}
&z_{14}\cdot z_{41}z_{42}\tss w=\frac{(n-1)(m-n)}{2n}z_{42}\tss w,
&& z_{24}\cdot z_{41}z_{42}\tss w=\frac{(n-1)(m+n+2)}{2n} z_{41}\tss w, \el\\
&z_{13}\cdot z_{41}z_{42}\tss w=0, && z_{23}\cdot z_{41}z_{42}\tss w=0.
\el
\end{align}

\section{Bases of $\kgot^+$-invariants in Kac modules}
\label{sec:abk}
\setcounter{equation}{0}

\noindent
$K^{+}(0,0)$:
\begin{align}
[0,0] \qquad
& w,\quad z_{41}z_{31}z_{42}z_{32}\tss w, &
\el \\
[1,1] \qquad
& z_{32}\tss w,\quad z_{31}z_{42}z_{32}\tss w, &
\el\\
[2,0] \qquad
& z_{42}z_{32}\tss w, &
\el\\
[0,2] \qquad
& z_{31}z_{32}\tss w.  &
\el
\end{align}

\noindent
$K^{+}(1,0)$:
\begin{align}
[1,0] \qquad
& w,\quad z_{41}z_{32}\tss w,\quad
z_{41}z_{31}z_{42}z_{32}\tss w,\el\\
[2,1] \qquad
& z_{32}\tss w,\quad z_{31}z_{42}z_{32}\tss w,\el\\
[0,1] \qquad
& z_{31}\tss w,\quad z_{41}z_{31}z_{32}\tss w,\el\\
[3,0] \qquad
& z_{42}z_{32}\tss w,  \el\\
[1,2] \qquad
& z_{31}z_{32}\tss w.
\el
\end{align}
\noindent
$K^{+}(0,1)$:
\begin{align} 
[0,1] \qquad
& w,\quad z_{41}z_{32}\tss w,\quad
z_{41}z_{31}z_{42}z_{32}\tss w, \el \\
[1,2] \qquad
& z_{32}\tss w,\quad z_{31}z_{42}z_{32}\tss w,\el\\
[1,0] \qquad
& z_{42}\tss w,\quad z_{41}z_{42}z_{32}\tss w,\el\\
[2,1] \qquad
& z_{42}z_{32}\tss w,  \el\\
[0,3] \qquad
& z_{31}z_{32}\tss w.
\el
\end{align}

\noindent
$K^{+}(1,1)$:
\begin{align}
[1,1] \qquad & w,\quad  z_{41}z_{31}z_{42}z_{32}\tss w,\quad
 z_{31}z_{42}\tss w,\quad z_{41}z_{32}\tss w, \el\\
[2,2] \qquad
& z_{32}\tss w,\quad z_{31}z_{42}z_{32}\tss w, \el\\
[2,0] \qquad
& z_{42}\tss w,\quad z_{41}z_{42}z_{32}\tss w, \el\\
[0,2] \qquad
& z_{31}\tss w,\quad z_{41}z_{31}z_{32}\tss w, \el\\
[0,0] \qquad
& z_{41}\tss w,\quad z_{41}z_{31}z_{42}\tss w, \el\\
[3,1] \qquad
& z_{42}z_{32}\tss w, \el\\
[1,3] \qquad
& z_{31}z_{32}\tss w. \el
\end{align}

\noindent
$K^{+}(m,0)$ with $m\geqslant 2$:
\begin{align}
[m,0] \qquad
& w,\quad z_{41}z_{32}\tss w,\quad
z_{41}z_{31}z_{42}z_{32}\tss w, \el\\
[m+1,1] \qquad
& z_{32}\tss w,\quad z_{31}z_{42}z_{32}\tss w, \el\\
[m-1,1] \qquad
& z_{31}\tss w,\quad z_{41}z_{31}z_{32}\tss w,\el\\
[m+2,0] \qquad
& z_{42}z_{32}\tss w, \el\\
[m,2] \qquad
& z_{31}z_{32}\tss w, \el\\
[m-2,0] \qquad
& z_{41}z_{31}\tss w.\el
\end{align}
\noindent
$K^{+}(0,n)$ with $n\geqslant 2$:
\begin{align}
[0,n] \qquad
& w,\quad z_{41}z_{32}\tss w,\quad
z_{41}z_{31}z_{42}z_{32}\tss w,  \el \\
[1,n+1] \qquad
& z_{32}\tss w,\quad z_{31}z_{42}z_{32}\tss w, \el\\
[1,n-1] \qquad
& z_{42}\tss w,\quad z_{41}z_{42}z_{32}\tss w, \el\\
[2,n] \qquad
& z_{42}z_{32}\tss w, \el\\
[0,n+2] \qquad
& z_{31}z_{32}\tss w, \el\\
[0,n-2] \qquad
& z_{41}z_{42}\tss w. \el
\end{align}

\noindent
$K^{+}(m,1)$ with $m\geqslant 2$:
\begin{align} 
[m,1] \qquad
& w,\quad z_{41}z_{32}\tss w,\quad  z_{31}z_{42}\tss w,\quad
z_{41}z_{31}z_{42}z_{32}\tss w, \el\\
[m+1,2] \qquad
& z_{32}\tss w,\quad z_{31}z_{42}z_{32}\tss w, \el\\
[m+1,0] \qquad
& z_{42}\tss w,\quad z_{41}z_{42}z_{32}\tss w, \el\\
[m-1,2] \qquad
& z_{31}\tss w,\quad z_{41}z_{31}z_{32}\tss w,\el\\
[m-1,0] \qquad
& z_{41}\tss w,\quad z_{41}z_{31}z_{42}\tss w,\el\\
[m+2,1] \qquad
& z_{42}z_{32}\tss w, \el\\
[m,3] \qquad
& z_{31}z_{32}\tss w, \el\\
[m-2,1] \qquad
& z_{41}z_{31}\tss w.
\el
\end{align}
\noindent
$K^{+}(1,n)$ with $n\geqslant 2$:
\begin{align}
[1,n] \qquad
& w,\quad z_{41}z_{32}\tss w,\quad  z_{31}z_{42}\tss w,\quad
z_{41}z_{31}z_{42}z_{32}\tss w, \el\\
[2,n+1] \qquad
& z_{32}\tss w,\quad z_{31}z_{42}z_{32}\tss w, \el\\
[2,n-1] \qquad
& z_{42}\tss w,\quad z_{41}z_{42}z_{32}\tss w, \el\\
[0,n+1] \qquad
& z_{31}\tss w,\quad z_{41}z_{31}z_{32}\tss w,\el\\
[0,n-1] \qquad
& z_{41}\tss w,\quad z_{41}z_{31}z_{42}\tss w,\el\\
[3,n] \qquad
& z_{42}z_{32}\tss w, \el\\
[1,n+2] \qquad
& z_{31}z_{32}\tss w, \el\\
[1,n-2] \qquad
& z_{41}z_{42}\tss w.
\el
\end{align}
\noindent
$K^{+}(m,n)$ with $m,n\geqslant 2$:
\begin{align}
[m,n] \qquad
& w,\quad z_{41}z_{32}\tss w,\quad  z_{31}z_{42}\tss w,\quad
z_{41}z_{31}z_{42}z_{32}\tss w, \el\\
[m+1,n+1] \qquad
& z_{32}\tss w,\quad z_{31}z_{42}z_{32}\tss w, \el\\
[m+1,n-1] \qquad
& z_{42}\tss w,\quad z_{41}z_{42}z_{32}\tss w, \el\\
[m-1,n+1] \qquad
& z_{31}\tss w,\quad z_{41}z_{31}z_{32}\tss w,\el\\
[m-1,n-1] \qquad
& z_{41}\tss w,\quad z_{41}z_{31}z_{42}\tss w,\el\\
[m+2,n] \qquad
& z_{42}z_{32}\tss w, \el\\
[m,n+2] \qquad
& z_{31}z_{32}\tss w, \el\\
[m-2,n] \qquad
& z_{41}z_{31}\tss w, \el\\
[m,n-2] \qquad
& z_{41}z_{42}\tss w. \el
\end{align}

\end{document}